\documentstyle[11pt]{article}
\def\tF{{\tilde F}}
\def\tG{{\tilde G}}
\def\tE{{\tilde E}}
\def\hF{{\hat F}}
\def\hG{{\hat G}}
\def\hE{{\hat E}}

\def\oF{{\overline{F}}}
\def\oG{{\overline{G}}}
\def\CC{{\cal C}}
\def\uF{{\underline{F}}}
\def\uG{{\underline{G}}}
\def\oC{{\overline{C}}}
\def\uC{{\underline{C}}}

\begin{document}




\title{On the adjustment coefficient, drawdowns and Lundberg-type bounds for random walk}

\setcounter{page}{0}

\baselineskip= 20pt

\author{
Isaac Meilijson\footnotemark \\
{\em School of Mathematical Sciences} \\
{\em Raymond and Beverly Sackler Faculty of Exact Sciences} \\
{\em Tel-Aviv University, 69978 Tel-Aviv, Israel} \\
{\em E-mail: \tt{MEILIJSON@MATH.TAU.AC.IL}} \\
}

\maketitle


\begin{abstract}
Consider a random walk whose (light-tailed) increments have positive
mean. Lower and upper bounds are provided for the expected maximal
value of the random walk until it experiences a given drawdown $d$.
These bounds, related to the Calmar ratio in Finance, are of the
form $(\exp\{\alpha d\}-1)/\alpha$ and $(K \exp\{\alpha
d\}-1)/\alpha$ for some $K>1$, in terms of the adjustment
coefficient $\alpha$ ($E[\exp\{-\alpha X\}]=1$) of the insurance
risk literature.
Its inverse ${1 \over
\alpha}$ has been recently derived by Aumann \& Serrano
as an index of riskiness of the random variable $X$.

This article also complements the Lundberg exponential stochastic
upper bound and the Cr\'{a}mer-Lundberg approximation for the
expected minimum of the random walk, with an exponential stochastic
lower bound. The tail probability bounds are of the form $C
\exp\{-\alpha x\}$ and $\exp\{-\alpha x\}$ respectively, for some
${1 \over K} < C < 1$.

Our treatment of the problem involves Skorokhod embeddings of random
walks in Martingales, especially via the Az\'{e}ma--Yor and Dubins
stopping times, adapted from standard Brownian Motion to exponential
Martingales.

\medskip


\end{abstract}

\righthyphenmin=55
\def\tF{{\tilde F}}
\def\tG{{\tilde G}}
\def\tE{{\tilde E}}
\def\hF{{\hat F}}
\def\hG{{\hat G}}
\def\hE{{\hat E}}
\def\oF{{\overline{F}}}
\def\oG{{\overline{G}}}
\def\CC{{\cal C}}
\def\uF{{\underline{F}}}
\def\uG{{\underline{G}}}
\def\oC{{\overline{C}}}
\def\uC{{\underline{C}}}
%
%
%




\section{Introduction}

{\bf Drawdowns of Brownian Motion with positive drift.} Let
$\{W(t) \mid t \ge 0,\, W(0) = 0\}$ be Standard Brownian Motion
(SBM) and let $\{B(t) \mid B(t) = \mu t + \sigma W(t),\, t \ge 0
\}$ be Brownian Motion (BM) with {\em drift} $\mu>0$ and {\em
diffusion parameter} $\sigma \in (0,\infty)$. For $d > 0$, define
the stopping time
\begin{equation} \label{stime1}
    \tau_d^{BM} =
    \min \{t | \max_{0 \le s \le t} B(s) \ge B(t)+d \}
\end{equation}
to be the first time to achieve a {\em drawdown}\/ of size $d$.
That is, $\tau_d^{BM}$ is the first time that BM has gone down by
$d$ from its record high value so far. As motivated by Taylor
\cite{Taylor}, an investor that owns a share whose value at time
$t$ is
$V_t = V_0 \exp(B(t))$, may consider selling it at time
$\tau_d^{BM}$ (for some $d > 0$) because it has lost for the first
time some fixed fraction $1-\exp(-d)$ of its previously held
highest value $V_0 \exp(M_d)$ (where $M_d=M_d^{BM} = \max_{0 \le s
\le \tau_d^{BM}} B(s) = B(\tau_d^{BM}) + d$), a possible
indication of change of drift.

As pointed out in Meilijson \cite{IMEdraw}, drawdowns are {\em
gaps} for Dubins \& Schwarz \cite{LED-S}), {\em extents} for
Goldhirsch \& Noskovicz \cite{Go-No}) and {\em downfalls} for
Douady, Shiryaev \& Yor \cite{DShY}). Taylor \cite{Taylor} (see
also \cite{IMEdraw}) presents a closed form formula for the joint
moment generating function of $\tau_d^{BM}$ and $B(\tau_d^{BM})$,
from which it follows that $M_d^{BM}$ is exponentially
distributed, with expectation
\begin{equation} \label{expm1}
E[M_d^{BM}] = {\sigma^2 \over {2 \mu}} (\exp\{{{2 \mu} \over
\sigma^2} d\} - 1).
\end{equation}

\noindent {\bf Maximum of Brownian Motion with negative drift.}
The maximum \linebreak $\max(BM)=\inf_{t>0} \{B(t)\}$ is well
known to have the exponential distribution
\begin{equation} \label{maxBM}
P(\max(BM)>x)=1 \wedge \exp\{-{{2 |\mu|} \over \sigma^2} x\}.
\end{equation}

This article contributes to the generalization of (\ref{expm1}) and
(\ref{maxBM}) from BM to random walks (RW). There is a rather vast
literature on the maximum of RW with negative drift. Kingman
\cite{Kingman} showed that $P(\max(RW)>x) \approx 1 \wedge
\exp\{-{{2 |\mu|} \over \sigma^2} x\}$ for small $\mu$, Siegmund
\cite{Siegmund} studied first order corrections to this
approximation via renewal-type overflow distributions and Chang \&
Peres \cite{ChangPeres} developed asymptotic expansions of
$P(\max(RW)>x)$ for the Gaussian case. Blanchet \& Glynn
\cite{BlGl1} improved on these approximations. In the insurance risk
literature, exponential bounds and approximations of $P(\max(RW)>x)$
are referred to as Lundberg's inequality or Cr\'{a}mer-Lundberg
approximations (see Asmussen's comprehensive treatise
\cite{Asmussen}).

This paper is methodologically different from the above; instead of
relying on change of measure and renewal theory, our setup involves
exponential martingales and Skorokhod embeddings, in a way
reminiscent of Wald's \cite{Wald} method for deriving the OC
characteristic of the Sequential Probability Ratio Test. As part of
the change, we will give up on trying to save the inaccurate role of
${{2 |\mu|} \over \sigma^2}$ as the exponential rate in the
questions under study, in favor of the so-called {\em adjustment
coefficient} of the insurance risk literature, provided by the
$\alpha$ solving $E[\exp\{-\alpha X\}]=1$. However, the rate ${{2
|\mu|} \over \sigma^2}$ will stay around: the RW will be coupled
with a BM for which ${{2 |\mu|} \over \sigma^2}$ is $\alpha$.

More explicitely, using Skorokhod (\cite{SKOR} and also
\cite{Az-Yo,Ch-Wa,LED,IME}) embeddings, mean-zero RW can be viewed
as optional sampling of SBM. This idea will be mimicked here to
embed the exponential Martingale $\exp\{-\alpha S_n\}$ into the
Martingale $\exp\{-\alpha B(t)\}$. This method could be useful in
obtaining other approximate extensions of pricing under log-normal
models to more general distributions.

Aumann \& Serrano \cite{AumSerr} asked a scalar index of {\em
riskiness} $Q(X)$ of the random variable (r.v.) $X$ to satisfy an
{\em homogeneity} axiom $Q(tX)=tQ(X)$ and a {\em duality} axiom that
models the increased preference of a more risk averse individual for
constant wealth $w$ over random wealth $w+X$. The unique solution
(up to a multiplicative constant) is the inverse ${1 \over \alpha}$
of the adjustment coefficient. The role played by $\alpha$ in our
subject matter is clearly consistent with riskiness - a large
$\alpha$ corresponds to low risk, as it ({\em i}) protects against
heavy initial losses before eventual divergence of the RW to
$\infty$, and ({\em ii}) makes the RW reach high yield before
experiencing sizable drawdowns.

The {\em Calmar ratio} (see Atiya \& Magdon-Ismail
\cite{AtiyaMagdon} and the implementation of their work in the
Matlab financial toolbox) of a financial asset with positive drift
is a measure of the likely drawdown in the logarithm of its price in
a given interval of time, such as a year. Since height (and time)
are exponential in the drawdown, the Calmar ratio is heavily
influenced by the length of this time interval. Besides, typical
drawdown in a given time span is harder to analyze than our subject
matter, typical height (or time) to achieve a given drawdown. We
propose the use of the adjustment coefficient or its inverse as a
Calmar-type measure of the risk of a financial asset, and provide
simple approximate formulas to quantify its effects. The more
commonly used {\em Sharpe index}, or ratio of net drift (drift minus
market interest rate) to volatility (standard deviation), lets
volatility penalize the asset even when it favors gains. In
contrast, drawdown-based indices measure risk in a more reasonable
asymmetric sense.

\section{Results}

From now on, we only consider BM and RW with positive drift and thus
unify the presentation of the two problems, by switching from the
commonly studied maximum of BM and RW with negative drift to the
equivalent treatment of the minimum of BM and RW with positive
drift.

Consider a r. v. $X \sim F$ with $0 < E[X] < \infty$ and $P(X<0)>0$.
Assume further the existence of $\alpha>0$ such that
$E[\exp\{-\alpha X\}]=1$. Since (if finite) the moment generating
function $\Psi(t)=E[\exp\{t X\}]$ is strictly convex with
$\Psi'(0)=E[X]>0$ and $\Psi(t) \rightarrow \infty$ as $|t|
\rightarrow \infty$, such $\alpha$ exists and is unique as long as
the moment generating function is finite wherever relevant. This
assumption is satisfied e.g. for Gaussian r.v.'s and for r.v.'s
bounded from below. If $X \sim N(\mu, \sigma^2)$, then $\alpha$ is
indeed ${{2 \mu} \over \sigma^2}$ (see (\ref{expm1})).

Besides $\alpha$, we need other characteristics of the
distribution $F$.
\begin{eqnarray} \label{d0}
d^{+} & = & {1 \over \alpha} \sup_{0<x<\mbox{es}_F} -
\log(E[e^{-\alpha(X-x)} | X \ge x]) \nonumber \\
d^{-} & = & {1 \over \alpha} \sup_{\mbox{ei}_F<x<0}
\log(E[e^{\alpha(x-X)} | X < x]) \\
d_0 & = & d^{+} + d^{-} \nonumber
\end{eqnarray}
where $\mbox{es}_F = \sup\{y|F(y)<1\}$ and $\mbox{ei}_F =
\inf\{y|F(y)>0\}$ are the {\em essential supremum} and {\em infimum}
of $F$. By Jensen's inequality, $d^{+}$ is bounded from above by the
simpler and more natural $\sup_x E[X-x | X \ge x]$ and $d^{-}$ is
accordingly bounded from below. These constants are defined in terms
of excesses of the r. v. $X$ itself, unlike the Siegmund or
Cr\'{a}mer-Lundberg approximations, built in terms of the renewal
overflow distribution of the random walk with $X$-increments.

\bigskip

Let $X_i \ ; \ i=1,2,\dots$ be i.i.d. $F-$distributed random
variables and let $S_0=0 \ ; \ S_n=\sum_{i=1}^n X_i$ be the
corresponding random walk (RW). The definition of $\alpha$ makes
$\exp\{-\alpha S_n\}$ a mean-1 Martingale.

Drawdowns, maximal heights $M_d^{RW}$ achieved prior to drawdowns
and the corresponding stopping times $\tau_d^{RW}$ can be defined
for random walk in much the same way they are defined for Brownian
Motion. Re-stating Lundberg's inequality as the RHS of
(\ref{distrmin}), the purpose of this paper is to prove the other
three inequalities in (\ref{expectm1}) and (\ref{distrmin})

\bigskip

\begin{equation} \label{expectm1}
{{e^{\alpha d} - 1} \over \alpha} \le E[M_d^{RW}] \le {{e^{\alpha
(d+d_0)} - 1} \over \alpha}
\end{equation}

For $x > 0$,

\begin{equation} \label{distrmin}
e^{- \alpha (x+d^{-})} \le P(-\min(RW)>x) \le e^{- \alpha x}
\end{equation}

We thus have lower and upper bounds for $E[M_d^{RW}]$ whose ratio
stays bounded as $d$ increases, provided $d_0$ is finite. These
bounds clearly show that drawdowns are logarithmic in the highest
value achieved so far, and precisely identify the exponential rate
$\alpha$ at which the latter grows as a function of the former.
Besides claiming bounds on the mean, the upper bound in
(\ref{expectm1}) is in fact a stochastic inequality: $M_d^{RW}$ is
stochastically smaller than the exponentially distributed random
variable $M_{d+d_0}^{BM}$. We do not provide a stochastic lower
bound. In contrast, (\ref{distrmin}) provides stochastic upper and
lower bounds on the minimum of RW.

The upper bound in (\ref{expectm1}) can be improved by letting $d_0$
depend on $d$ and be defined as $d_0$ in (\ref{d0}) but restricting
the maximization to \linebreak
$x \in (0, \min(es_F,d))$ and $x \in
(\max(ei_F,-d),0)$. This improved upper bound is finite for every
$d$.

\section{A few examples}

{\large {\bf Example 1: The Gaussian case.}} Let X have a normal
distribution with positive mean $\mu$ and standard deviation
$\sigma$. Let $\phi$ and $\Phi$ stand respectively for the
standard normal density and cumulative distribution function. Then
\begin{equation} \label{egnormal}
\alpha = {{2 \mu} \over \sigma^2} \ ; \ e^{\alpha d^{+}} =
e^{\alpha d^{-}} ={{\Phi({\mu \over \sigma})} \over {1-\Phi({\mu
\over \sigma})}}
\end{equation}
and this proves the two following rather elegant formulas.
\begin{eqnarray}
{{e^{{{2 \mu} \over \sigma^2} d}-1} \over {{2 \mu} \over \sigma^2}}
& \le & \hspace{1.1cm} E[M_d] \hspace{1.2cm} \le {({{\Phi({\mu \over
\sigma})} \over {1-\Phi({\mu \over \sigma})}})^2 {e^{{{2 \mu} \over
\sigma^2} d}-1} \over {{2
\mu} \over \sigma^2}} \label{Gaussdraw} \\
{{1-\Phi({\mu \over \sigma})} \over \Phi({\mu \over \sigma})}
e^{-{{2 \mu} \over \sigma^2} x} & \le & P(-\min(RW) > x) \le e^{-{{2
\mu} \over \sigma^2} x} \label{Gaussmin}
\end{eqnarray}

If we view the normal random walk as sampling Brownian Motion with
drift $\mu$ and diffusion coefficient $\sigma$ at regular intervals,
$\alpha$ is independent of the grid length $\delta$ but the three
$d$'s are not, predictably vanishing with the grid length.

The LHS of (\ref{egnormal}) is well known and easy to obtain from
the formula $\exp\{\mu+\sigma^2 t^2/2\}$ of the moment generating
function of the normal distribution. As for the RHS, it requires
evaluating via
\begin{equation} \label{normalplusZ}
E[e^{-\beta Z} | Z > z] = e^{{1 \over 2} \beta^2} {1 \over \sqrt{2
\pi}} {{\int_z^\infty \exp\{-{1 \over 2} (z+\beta)^2\} dz} \over
{1-\Phi(z)}} = e^{{1 \over 2} \beta^2} {{1-\Phi(z+\beta)} \over
{1-\Phi(z)}}
\end{equation}
the expressions
\begin{eqnarray}
E[e^{-\alpha X} | X > x] & = & E[e^{-\alpha (\mu+\sigma Z)} |
\mu+\sigma Z > x] = e^{-2{{\mu^2} \over {\sigma^2}}} E[e^{-2 { \mu
\over \sigma} Z} | Z > {{x-\mu} \over \sigma}] \nonumber \\
& = & {{1-\Phi({{x+\mu} \over \sigma})} \over {1-\Phi({{x-\mu}
\over \sigma})}} \label{normalplus} \\
E[e^{-\alpha X} | X < x] & = & {{\Phi({{x+\mu} \over \sigma})}
\over {\Phi({{x-\mu} \over \sigma})}} \label{normalminus}
\end{eqnarray}
from which the RHS of (\ref{egnormal}) follows, at least in the
sense of plugging $x=0$. To see that $x=0$ is indeed the correct
choice for each side, observe that the normal distribution is IFR
- has increasing failure rate (Mills' ratio ${\phi(z) \over
{1-\Phi(z)}}$). But IFR implies that the residual distributions
${\cal L}(X-x|X>x)$ are ordered by stochastic inequality. Hence,
so are the expectations of monotone functions, such as the
exponential function. This argument applies equally to the two
tails.

\bigskip

\noindent {\large {\bf Example 2: The double exponential case.}}
Let $X$ have density $p \theta \exp(-\theta x\}$ for $x>0$ and
$(1-p) \mu \exp(\mu x\}$ for $x<0$, with ${\theta \over
{\mu+\theta}} < p < 1$. Then
\begin{equation} \label{doublexp1}
E[X]={p \over \theta}-{{1-p} \over \mu} \ ; \ \alpha = p \mu -
(1-p) \theta
\end{equation} with the corresponding bound ingredients
\begin{equation} \label{doublexp2}
e^{\alpha d^{+}} = {{(\mu+\theta) p} \over \theta} \ ; \ e^{\alpha
d^{-}} = {\mu \over {(\mu+\theta) (1-p)}} \ ; \ e^{\alpha d_0} =
{{p \mu} \over {(1-p) \theta}}
\end{equation}

The rate $\alpha$ exceeds ${{2 E[X]} \over \mbox{Var}[X]}$ for all
$p$ if $\mu \ge \theta$, but if $\mu < \theta$ the opposite
inequality holds for all $p$ close enough to $1$.

\bigskip

\noindent {\large {\bf Example 3: The shifted exponential case.}}
Let $X$ have exponential distribution with mean ${1 \over \theta}$
shifted down by $\Delta < {1 \over \theta}$ so as to allow negative
values and still preserve positive mean. It is easier to express the
inverse function to $\alpha$:
\begin{equation} \label{egexponential1}
{1 \over \alpha} \log(1+{\alpha \over \theta}) = \Delta
\end{equation}
from which
\begin{equation} \label{egexponential2}
d^{+} = \Delta \ ; \ e^{\alpha d^{+}} = 1 + {\alpha \over \theta}
\ ; \ e^{\alpha d^{-}} = {{1-e^{-(\theta+\alpha) \Delta}} \over
{1-e^{-\theta \Delta}}} \ ; \ e^{\alpha d_0} = 1 + {{\alpha \over
\theta} \over {1-e^{-\theta \Delta}}}
\end{equation}

The Gaussian-motivated rate ${{2 E[X]} \over \mbox{Var}[X]}$ is $2
\theta (1- \theta \Delta)$, always smaller than $\alpha$. That is,
a random walk with shifted exponential increments gets to higher
heights before a given drawdown than a normal one with the same
mean and variance.

\bigskip

\noindent {\large {\bf Example 4: A dichotomous case.}} Let
$P(X=-1)=1-p$ and \linebreak
$P(X=1)=p > {1 \over 2}$. Then $\alpha =
\log {p \over {1-p}}$ and, obviously, $d^{+}=d^{-}=1$. As is well
known from the Gambler's ruin problem, the probability of reaching
$+1$ before (integer) $-d$ is ${{1 - \exp\{-\alpha d\}} \over
{1-\exp\{-\alpha (d+1)\}}}$. $M_d^{RW}$ is nothing but the number of
independent such attempts until a first ``failure''. Hence, it is
($-1$ plus) a geometric r. v., and its mean is
\begin{equation} \label{egdichot}
E(M_d^{RW}) = -1 + {1 \over \mbox{probability}} = {p \over {2 p
-1}} (e^{\alpha d}-1) \ .
\end{equation}
For non-integer $d$, the ceiling of $d$ should be substituted in
(\ref{egdichot}). Even without doing so, the LHS of
(\ref{expectm1}) is verified, because $\alpha = \log{p \over
{1-p}}
> {{2 p -1} \over p}$. To ascertain the RHS, take $E(M_{d+1}^{RW})$ as
worst-case ceiling and check that ${p \over {2 p -1}} (e^{\alpha
(d+1)}-1)$ is below the bound $(\exp\{\alpha (d+2)\}-1)/\alpha$.

Just as in the shifted exponential case, the rate of growth
$\alpha$ exceeds the rate ${{2 E[X]} \over \mbox{Var}[X]} = {1
\over 2}({1 \over {1-p}} - {1 \over p})$ that would have obtained
in the Gaussian case. However,

\bigskip

\noindent {\large {\bf Example 5: A skew dichotomous case.}} Let
$P(X=-1)={b \over {1+b}}$ and \linebreak $P(X=b(1+\epsilon))={1
\over {1+b}}$. The mean is $E[X]={{\epsilon b} \over {1+b}}$ so let
us take \linebreak
$\epsilon=0.2$ to achieve positive mean and
$b=0.1$ to tilt the distribution towards bigger losses. Then
$\alpha=0.318$ but ${{2 E[X]} \over {\mbox{Var}[X]}} = 0.351$. This
shows that even for dichotomous variables the inequality between the
two can go both ways. Examples 2 and 5 suggest that
yield-to-drawdown performance worse than Gaussian is obtained when
the left tail is heavier than the right tail.

\bigskip

In all the previous examples, the distribution $F$ has
non-decreasing failure rate and the ``excess lifetime'' over $x$
looks shorter as $x$ increases. That's why the $d$'s are attained
at $x=0$ (see (\ref{d0})). It is easy to produce a four-point
distribution with one negative atom in which $d^{+}$ will be the
distance between the two rightmost atoms.

\bigskip

\noindent {\large {\bf Example 6: A power-law right tail.}} If $F$
is light left tailed but behaves like power law at the right tail,
then $\alpha$ is finite but $d^{+}$ is infinite because its
maximand behaves like $\log x$. To wit,
\begin{equation} \label{powerlaw}
E[e^{-\alpha (X-x)}|X>x]= {\gamma \over x} \int_0^\infty
{e^{-\alpha t} \over (1+{t \over x})^{\gamma+1}} dt \approx
{\gamma \over {\alpha x}}
\end{equation}
so $-{1 \over \alpha} \log(E[\exp\{-\alpha (X-x)\}|X>x]) = \log(x)
+ \mbox{o}(1)$. Although much smaller than
$E[X-x|X>x]=\mbox{O}(x)$ (see the sentence following (\ref{d0})
and the remark in the next section), it still goes to $\infty$.
However, the improved definition of $d^{+}$ sets it as $\log(d)$
up to a vanishing term.

This example illustrates that yield-to-drawdown, while at least as
high as the Brownian lower bound, may in principle be
superexponential.

\section{Miscellaneous}

{\bf The record high value $M_d^{BM}$ is exponentially
distributed}. This is so because as long as first hitting times of
positive heights occur before achieving a drawdown of $d$, these
times are renewal times: knowing that $M_d^{BM} > x$ is the same
as knowing that $B$ has not achieved a drawdown of $d$ by the time
it first reached height $x$. But then it starts anew the quest for
a drawdown.

\bigskip

\noindent {\bf A direct argument for (\ref{expm1})}. Since the
mean-$1$ Martingale $\exp\{-\alpha B\}$ stopped at $\tau_d^{BM}$
is uniformly bounded, it is also uniformly integrable. Hence,
\begin{equation} \label{direct}
1=E[e^{-\alpha B(\tau_d^{BM})}]= E[e^{-\alpha M_d^{BM}}] e^{\alpha
d}
\end{equation}
Since $M_d^{BM}$ is exponentially
distributed, $E[\exp\{-\alpha M_d^{BM}\}]= {1 \over {1+ \alpha
E[M_d^{BM}]}}$.

\section{Skorokhod embeddings in Martingales}

The problem as posed and solved by Skorokhod in \cite{SKOR} is the
following: given a distribution $F$ of a r. v. $Y$ with mean zero
and finite variance, find a stopping time $\tau$ in SBM $W$, with
finite mean, for which $W(\tau)$ is distributed $F$. The
Chacon--Walsh \cite{Ch-Wa} family of solutions is easiest to
describe: Express $Y$ as the limit of a Martingale
$Y_n=E[Y|{\cal{F}}_n]$ with dichotomous transitions (that is, the
conditional distribution of $Y_{n+1}$ given ${\cal{F}}_n$ is a.s.
two-valued), and then
\linebreak
progressively embed this Martingale
in $W$ by a sequence of first exit times from open intervals.

\bigskip

Dubins \cite{LED} was the first to build such a scheme, letting
${\cal F}_1$ decide whether $Y \ge E[Y]$ or $Y < E[Y]$ by a first
exit time of $W$ starting at $E[Y]$ from the open interval $(E[Y|
Y < E[Y]], E[Y|Y \ge E[Y]])$. It then proceeds recursively. E.g.,
if the first step ended at $E[Y|Y \ge E[Y]]$ then the second step
ends when $W$, re-starting at $E[Y|Y \ge E[Y]]$, first exits the
open interval $(E[Y|E[Y] \le Y < E[Y|Y \ge E[Y]]], E[Y|Y \ge E[Y|
Y \ge E[Y]]])$.

\bigskip

One of the analytically most elegant solutions to Skorokhod's
problem is the Az\'{e}ma--Yor stopping time $T_{AY}$ (see Az\'{e}ma
\& Yor \cite{Az-Yo} and Meilijson \cite{IME}), defined in terms of
$H_F(x)=E[Y|Y \ge x]=\int_x^\infty y dF(y)/(1-F(x-))$, the upper
barycenter function of $F$, as
\begin{equation} \label{AYstop}
T_{AY}=\min\Bigl\{t \Bigm| \max_{0 \le s \le t} W(s) \ge
H_F\bigl(W(t)\bigr)\Bigr\}.
\end{equation}

Among all uniformly integrable c\`{a}dl\`{a}g Martingales with a
given final or limiting distribution, SBM stopped at the
Az\'{e}ma--Yor stopping time to embed this distribution is extremal,
in the sense that it stochastically maximizes the maximum of the
Martingale (see Dubins \& Gilat \cite{LED-DG} and Az\'{e}ma \& Yor
\cite{Az-Yo}). That is, if $T_{AY}$ embeds $F$ then $M_{T_{AY}}$ is
stochastically bigger than the maximum of any such Martingale.

The connection of the Az\'{e}ma--Yor stopping time to the
Chacon--Walsh family becomes apparent (see Meilijson \cite{IME}) if
the r. v. $Y$ has finite support $\{x_1 < \cdots < x_k\}$. In this
case, let ${\cal{F}}_n$ be the $\sigma$-field generated by
$\min(Y,x_{n+1})$, that is, let the atoms of $Y$ be incorporated one
at a time, in their natural order: the first stage decides whether
$Y=x_1$ (by stopping there) or otherwise (by temporarily stopping at
$E[Y|Y>x_1]$), etc. This is precisely the Az\'{e}ma--Yor stopping
rule: stop as soon as a value of $Y$ is reached after having visited
the conditional expectation of $Y$ from this value and up.

Clearly, there is a mirror-image notion $T_{AY-}$ to Az\'{e}ma \&
Yor's stopping time that stochastically minimizes the minimum of
the Martingale. Simply put, apply $T_{AY}$ to embed the
distribution of $-X$ in $-W$.

\bigskip

The stopping time $T_{DAY}$ to be applied here is a hybrid of the
Dubins and Az\'{e}ma \& Yor stopping times. It starts as the
Dubins stopping time by a first-exit time of SBM $W$ from the
interval $(E[Y| Y < E[Y]], E[Y|Y \ge E[Y]])$. If exit occurred at
the top it proceeds by embedding the law ${\cal L}(Y|Y \ge E[Y])$
by $T_{AY}$ in the remainder SBM starting at $E[Y|Y \ge E[Y]]$.
If, on the other hand, exit occurred at the bottom it proceeds by
embedding the law ${\cal L}(Y|Y < E[Y])$ by $T_{AY-}$ in the
remainder SBM starting at $E[Y|Y < E[Y]]$.

\bigskip

Once a distribution $F$ is embeddable in SBM W, so is the random
walk with increments distributed $F$. Plainly, embed $X_1$ at time
$\tau_1$, then use the same rule to embed $X_2$ at time $\tau_2$
in the SBM $W'(t)=W(\tau_1+t)-W(\tau_1)$, etc. Skorokhod's
original idea was to infer the Central Limit Theorem for ${S_n
\over \sqrt{n}}$ from the Law of Large Numbers for ${\sum_{i=1}^n
\tau_i \over n}$. This idea was extended by Holewijn \& Meilijson
\cite{HoleMei} from random walks to Martingales with stationary
ergodic increments, to obtain a simple proof of the Billingsley \&
Ibragimov \cite{BillIbr} CLT.

\bigskip

Just as a random walk $S_n$ can be embedded in SBM, the
exponential Martingale $\exp\{-\alpha S_n\}$ can be embedded in
the continuous-time continuous Martingale $\exp\{-\alpha B(t)\}$,
where the BM $B$ has drift $\mu$ and diffusion coefficient
$\sigma$ such that ${{2 \mu} \over \sigma^2} = \alpha$. At the
time the RW reaches drawdown at least $d$, BM has also gone down
by at least $d$, but may have gone higher in the meantime. Thus,
$\tau_d^{BM}$ is a.s. smaller than $\tau_d^{RW}$. Now we may
compute, under the obvious definition of $\Delta \ge d$ a.s.,
\begin{eqnarray} \label{BMandRW}
\mu E[\tau_d^{BM}] & = & E[B(\tau_d^{BM})] = E[M_d^{BM}] - d \nonumber \\
& \le & \mu E[\tau_d^{RW}] = E[B(\tau_d^{RW})] =
E[M_d^{RW}]-E[\Delta]
\end{eqnarray}
so $E[M_d^{RW}] \ge E[M_d^{BM}]$. We have proved the LHS
inequality in (\ref{expectm1}) by the method of Coupling.

The RHS inequality in (\ref{expectm1}) is proved by using the
Dubins - Az\'{e}ma \& Yor stopping time $\tau_{DAY}$ for the above
embedding. If this embedding in $\exp\{-\alpha W\}$ ends up with
$W$ below (resp. above) some $x > 0$ (non-positive), the
underlying process could not have reached the exponential
barycenter height above (below) the closest support point to the
right (left) of $x$, because then the stopped value would have
been from this support rightwards (leftwards). For the first
increment of RW following the maximal (minimal) value the relevant
$x$ is 0, but for values embedded later the starting $x$ is lower
(higher). It should now be clear that the BM path can't reach as
far up as $M_d^{RW} + d^{+}$ nor as far down as $B(\tau_d^{RW}) -
d^{-}$ before RW achieves drawdown $d$. Hence, $BM$ can't reach
drawdown $d+d_0$ before $RW$ reaches drawdown $d$, or
$\tau_{d+d_0}^{BM} \ge \tau_d^{RW}$ a.s. The RHS inequality in
(\ref{expectm1}) follows.

This RHS inequality holds stochastically, since the cumulative
maximum of $BM$ exceeds the cumulative maximum of the embedded
$RW$ timewise. Thus, it holds a fortiori if the former is measured
later than the latter. This argument fails for the LHS because
then the latter is measured before the former.

\bigskip

Finally, a non-trivial stochastic upper bound on $M_d^{RW}$ can be
defined even if $d_0=\infty$, as announced in the Introduction. In
Example 5, this corresponds to adding a little over $\log d$ to
the exponent, or, the upper bound is roughly $d$ times the lower
bound.

\bigskip

\section*{Acknowledgements}
I wish to thank Linir Or for discussions that led to this study and
David Gilat and David Siegmund for useful comments, for which I also
thank Natalia Lysenko, who introduced me to Calmar rates. The
hospitality of ETH Z\"{u}rich is warmly appreciated.

\end{document}